\documentclass[11pt]{article}
\newcommand{\qed}{{\hfill\rule{4pt}{7pt}}}
\usepackage{amsmath, epsfig}
\usepackage{amssymb}
\usepackage{amsfonts}
\usepackage{latexsym}
\newtheorem{thm}{Theorem}

\newtheorem{cor}[thm]{Corollary}

\def\pf{\noindent {\it Proof.} }

\begin{document}

\begin{center}
{\Large\bf  A note on two identities arising from enumeration of convex polyominoes}
\end{center}

\vskip 2mm
\centerline{Victor J. W. Guo$^1$ and Jiang Zeng$^{1,2}$}

\begin{center} $^1$ Center for Combinatorics, LPMC\\
Nankai University, Tianjin 300071, People's Republic of China\\
{\tt jwguo@eyou.com}
\end{center}

\begin{center}
$^2$ Institut Girard Desargues,
Universit\'e Claude Bernard (Lyon I)\\
F-69622, Villeurbanne Cedex, France \\
{\tt zeng@desargues.univ-lyon1.fr}
\end{center}

\vskip 0.7cm \noindent{\bf Abstract.} {\small Motivated by some
binomial coefficients identities encountered in our  approach to
the enumeration of convex polyominoes, we prove  some more general
identities of the same type, one of which turns out to be related
to a strange evaluation of ${}_3F_2$ of Gessel and Stanton. }



\section{Introduction}
In our elementary approach to the enumeration of convex
polyominoes with an $(m+1)\times (n+1)$ minimal bounding
rectangle~\cite{GZ}, we encountered the  following two interesting
identities:
\begin{align}
&\sum_{a=1}^m\sum_{b=1}^n{m+n-a+b-1\choose m-a}{m+n+a-b-1\choose
n-b}
=\frac{mn}{2(m+n)}{m+n\choose m}^2, \label{eq:S3}\\[5pt]
&\sum_{a=1}^{m-2}\sum_{b=1}^{n-2}{m+n+a-b-1\choose m+a+1}
{m+n-a+b-1\choose n+b+1} \nonumber\\[5pt]
&={m+n\choose m}^2+{m+n\choose m-1}{m+n\choose
n-1}+\frac{mn}{2(m+n)}{m+n\choose m}^2-{2m+2n\choose 2n}. \label{eq:S4}
\end{align}
Although the single-sum case of binomial coefficients identities
is well-studied, the symbolic manipulation of binomial
multiple-sum identities depends on the performance of computers
(see, for example, \cite{AP}). Therefore, formulas of binomial
double-sums are still a challenge both for human and computer.

In this paper, we will give some generalizations of the
above two formulas, our main results are the following two theorems.
\begin{thm}\label{thm:mnalpha}
For $m,n\in \mathbb{N}$ and any number $\alpha\neq 0$, we have
\begin{align}
&\sum_{a=1}^{m}\sum_{b=1}^{n}{(1+\alpha)m-a+b-1\choose m-a}
{(1+\alpha^{-1})n+a-b-1\choose n-b}\nonumber\\[5pt]
&=\frac{mn}{(1+\alpha)(m+\alpha^{-1}n)} {(1+\alpha)m\choose m}
{(1+\alpha^{-1})n\choose n}. \label{eq:mnalpha}
\end{align}
\end{thm}

\begin{thm}For $m,n,r\in \mathbb{N}$ and any number $\alpha\neq 0$, we have
\begin{align}
&\sum_{a=0}^{m-r-2}\sum_{b=0}^{n-r-2}{(1+\alpha)m-a+b-1\choose
m-r-2-a}
{(1+\alpha^{-1})n+a-b-1\choose n-r-2-b}\nonumber\\[5pt]
&\hspace{2cm}+\sum_{a=0}^{m+r}\sum_{b=0}^{n+r}{(1+\alpha)m-a+b-1\choose
m+r-a}{(1+\alpha^{-1})n+a-b-1\choose n+r-b}\nonumber\\[5pt]
&=\frac{2mn}{(1+\alpha)(m+\alpha^{-1}n)} {(1+\alpha)m\choose
m}{(1+\alpha^{-1})n\choose
n}\nonumber\\[5pt]
&\hspace{2cm}+\sum_{k=-r}^{r}(r-|k|+1){(1+\alpha)m\choose
m-k}{(1+\alpha^{-1})n\choose n-k}. \label{eq:genrr-2}
\end{align}
\end{thm}

We use essentially the generating function 
techniques, that is, to prove $A=B$ we show that their generating functions
are equal. Two proofs of Theorem 1 will be given in Section 1, and Theorem 2 will be
proved in Section 3. In Section 4, we derive some interesting special cases of Theorems 1
and 2.

\section{Proofs of Theorem~1}

\subsection{First Proof of Theorem \ref{thm:mnalpha}.} Multiplying the left-hand side of \eqref{eq:mnalpha} 
by $x^my^n$ and summing over $m\geq 0$ and $n\geq 0$ we obtain the generating function
$F(x,y)$, which can be written after exchanging the order of summations as:
\begin{align*}
F(x,y) :=\sum_{a,b=1}^{\infty}x^ay^b\sum_{m,n=0}^{\infty}
{(1+\alpha)m+a\alpha+b-1\choose m}
{(1+\alpha^{-1})n+a+b\alpha^{-1}-1\choose n}x^my^n.
\end{align*}
Summing the two inner sums by the following classical formula (see
\cite[p. 146]{PS}) and \cite[(9)]{Gould}):
\[
\sum_{n=0}^{\infty}{\alpha+\beta n\choose
n}w^n=\frac{z^{\alpha+1}}{(1-\beta)z+\beta},\quad
\hbox{where}\quad  w=\frac{z-1}{z^{\beta}},
\]
 and then summing the two resulted geometric series over $a$ and $b$ we obtain
\begin{align*}
F(x,y)=\frac{uv(u-1)(v-1)} {(1+\alpha-\alpha
u)(1+\alpha^{-1}-\alpha^{-1}v)(uv-u-v)^2},
\end{align*}
where
\begin{equation}
x=\frac{u-1}{u^{1+\alpha}},\qquad
y=\frac{v-1}{v^{1+\alpha^{-1}}}.\label{eq:xuyv}
\end{equation}
Now, using the fact that $\sum_{k\geq 0}kx^k=x/(1-x)^2$ we have
\begin{align*}
F(x,y)&=\frac{uv}
 {(1+\alpha-\alpha u)(1+\alpha^{-1}-\alpha^{-1}v)} \sum_{k=0}^{\infty}k(u-1)^{k}(v-1)^{k} \\[5pt]
&=\sum_{k=0}^{\infty}\sum_{m,n=k}^{\infty} k{(1+\alpha)m\choose
m-k}{(1+\alpha^{-1})n\choose n-k}x^m y^n.
\end{align*}
Comparing with \eqref{eq:mnalpha}, it remains to check the following identity:
\begin{align}
\sum_{k=0}^{\min\{m,n\}} k{(1+\alpha)m\choose m-k}{(1+\alpha^{-1})n\choose n-k}
=\frac{mn}{(1+\alpha)(m+\alpha^{-1}n)} {(1+\alpha)m\choose
m}{(1+\alpha^{-1})n\choose n}. \label{eq:kxyalpha}
\end{align}
Writing
\begin{align*}
k{(1+\alpha)m\choose m-k}{(1+\alpha^{-1})n\choose n-k}
&=\frac{(m+\alpha^{-1}k)(n+\alpha k)}{(1+\alpha)(m+\alpha^{-1}n)}
{(1+\alpha)m\choose m-k}{(1+\alpha^{-1})n\choose n-k}\\[5pt]
&\quad\
-\frac{(m+\alpha^{-1}(k+1))(n+\alpha(k+1))}{(1+\alpha)(m+\alpha^{-1}n)}
{(1+\alpha)m\choose m-k-1}{(1+\alpha^{-1})n\choose n-k-1},
\end{align*}
equation~\eqref{eq:kxyalpha} follows then by summing over $k$ from $0$ to $\min\{m,n\}$.
 \qed

{\noindent\bf Remark}: Notice that
\begin{align*}
\sum_{m=0}^{\infty}{(1+\alpha)m\choose
m}x^m=\frac{u}{1+\alpha-\alpha u},\qquad
\sum_{n=0}^{\infty}{(1+\alpha^{-1})n\choose n}x^n
=\frac{v}{1+\alpha^{-1}-\alpha^{-1} v}.
\end{align*}
Multiplying the two sides of \eqref{eq:mnalpha} by
$(1+\alpha)(m+\alpha^{-1}n)$, we see that \eqref{eq:mnalpha} is
equivalent to:
\begin{align}
&\hskip -2mm[(1+\alpha)x\frac{\partial }{\partial x}
+(1+\alpha^{-1})y\frac{\partial }{\partial y}]F(x,y)
=xy\frac{d}{dx}\left(\frac{u}{1+\alpha-\alpha u}\right)
\frac{d}{dy}\left(\frac{v}{1+\alpha^{-1}-\alpha^{-1} v}\right).
\label{eq:xyalpha}
\end{align}
It's then possible to give another proof of \eqref{eq:mnalpha} by
checking \eqref{eq:xyalpha}, which is left to the interested reader.

\subsection{ Second Proof of Theorem \ref{thm:mnalpha}.} Replacing
$b$ by $k+1$ and writing the $k$-sum in standard hypergeometric
notation we can write the left-hand side as
\begin{align*}
L&:=\sum_{a=1}^{m}\sum_{k\geq 0} {(1+\alpha)m-a+k\choose m-a}
{(1+\alpha^{-1})n+a-2-k\choose \alpha^{-1}n+a-1}\\[5pt]
&=\sum_{a=1}^{m}{(1+\alpha)m-a\choose m-a}{(1+\alpha^{-1})n+a-2\choose n-1}\\[5pt]
&\quad\cdot {}_3F_2\left[\begin{array}{c}
  1-n,\ 1,\ (1+\alpha)m+1-a\\[5pt]
  \alpha m+1,\ -(1+\alpha^{-1})n-a+2
\end{array};1\right].
\end{align*}
Applying the transformation~\cite[p.~142]{AAR}:
$$ {}_3F_2\left[\begin{array}{c}
  -N,\ a,\ b\\[5pt]
   d,\ e
\end{array};1\right]=\frac{(e-b)_N}{(e)_N}
{}_3F_2\left[\begin{array}{c}
  -N,\ b,\ d-a\\[5pt]
   d,\ 1+b-e-N
\end{array};1\right]
$$
to the above ${}_3F_2$ we get
\begin{align*}
L&=\sum_{a=1}^{m}{(1+\alpha)m-a\choose
m-a}{(1+\alpha^{-1})n+a-2\choose n-1}
\frac{(-(1+\alpha)(m+\alpha^{-1}n)+1)_{n-1}}{(-(1+\alpha^{-1})n-a+2)_{n-1}}\\[5pt]
&\quad\ \cdot{}_3F_2\left[\begin{array}{c}
  1-n,\ (1+\alpha)m+1-a,\ \alpha m\\[5pt]
   \alpha m+1,\ (1+\alpha)m+\alpha^{-1}n+1
\end{array};1\right].
\end{align*}
Expanding the ${}_3F_2$ as a $k$-sum  and exchanging the order
with $a$-sum yields
\begin{align}
&L={(1+\alpha)m\choose m}{(1+\alpha^{-1})n-2\choose n-1}
\frac{(-(1+\alpha)(m+\alpha^{-1}n)+1)_{n-1}}{(2-(1+\alpha^{-1})n)_{n-1}}\nonumber\\[5pt]
&\quad\ \cdot\sum_{k\geq 0}\frac{(1-n)_k(\alpha
m)_k((1+\alpha)m+1)_k}
{(\alpha m+1)_k((1+\alpha)m+\alpha^{-1}n+1)_k k!}\frac{m}{(1+\alpha)m+k}\\[5pt]
&\quad\ \cdot{}_2F_1\left[\begin{array}{c}
  1-m,\ 1 \\[5pt]
   1-(1+\alpha)m-k
\end{array};1\right]\nonumber\\[5pt]
&={(1+\alpha)m\choose m}{(1+\alpha^{-1})n-2\choose n-1}
\frac{(-(1+\alpha)(m+\alpha^{-1}n)+1)_{n-1}}{(2-(1+\alpha^{-1})n)_{n-1}}\nonumber\\[5pt]
&\quad\ \cdot \frac{m}{\alpha m+1}{}_3F_2\left[
\begin{array}{ccc}
  1-n,\ \alpha m,\ (1+\alpha)m+1 \\[5pt]
   \alpha m+2,\ (1+\alpha)m+\alpha^{-1}n+1
\end{array};1
\right].\label{eq:3F2sum}
\end{align}
The theorem then follows by applying Gessel and Stanton's
formula~\cite[(1.9)]{GS}:
$$
{}_3F_2\left[\begin{array}{c}
  -sb+s+1,\ b-1,\ -N \\[5pt]
  b+1,\ s(-N-b)-N
\end{array};1\right]=\frac{(1+s+sN)_Nb(N+1)}{(1+s(b+N))_N(b+N)},
$$
with $N=n-1$, $b=\alpha m+1$ and $s=-1-\alpha^{-1}$. \qed

\noindent{\it Remark.} If $\alpha=1$, we can also evaluate the
${}_3F_2$ in \eqref{eq:3F2sum} by applying Dixon's
formula~\cite[p.~143]{AAR}:
\begin{align*}
&{}_3F_2\left[\begin{array}{c}
  a,\ b,\ c \\[5pt]
  1+a-b,\ 1+a-c
\end{array};1\right]\\[5pt]
&=\frac{\Gamma(1+a-b)\Gamma(1+a-c)\Gamma(1+\frac{a}{2})\Gamma(1+\frac{a}{2}-b-c)}
{\Gamma(1+a)\Gamma(1+\frac{a}{2}-b)
\Gamma(1+\frac{a}{2}-c)\Gamma(1+a-b-c)},
\end{align*}
and if $m=n$, we can apply Whipple's formula~\cite[p.~149]{AAR}:
\begin{align*}
&{}_3F_2\left[\begin{array}{c}
  a,\ 1-a,\ c \\[5pt]
 d,\ 1+2c-d
\end{array};1\right]\\[5pt]
&=\frac{2^{1-2c}\pi\Gamma(d)\Gamma(1-2c+d)}
{\Gamma(\frac{1}{2}+\frac{a}{2}+c-\frac{d}{2})
\Gamma(\frac{a}{2}+\frac{d}{2})\Gamma(1-\frac{a}{2}+c-\frac{d}{2})
\Gamma(\frac{1}{2}-\frac{a}{2}+\frac{d}{2}) }.
\end{align*}
\section{Proof of Theorem~2}
Consider the generating function:
\begin{align*}
G_r(x,y) &:=\sum_{m,n=-r}^{\infty}\sum_{a=0}^{m+r}\sum_{b=0}^{n+r}
{(1+\alpha)m-a+b-1\choose m+r-a} {(1+\alpha^{-1})n+a-b-1\choose
n+r-b}x^my^n.
\end{align*}
Using \eqref{eq:xuyv}, as the first proof of Theorem
\ref{thm:mnalpha}, we have
\begin{align}
G_r(x,y)=\frac{uv(u-1)^{-r}(v-1)^{-r}} {(1+\alpha-\alpha
u)(1+\alpha^{-1}-\alpha^{-1}v)(uv-u-v)^2}. \label{eq:f1xyr}
\end{align}
Replacing $r$ by $-r-2$ in \eqref{eq:f1xyr}, we obtain
\begin{align*}
G_{-r-2}(x,y):&=\sum_{m,n=r+2}^{\infty}\sum_{a=0}^{m-r-2}\sum_{b=0}^{n-r-2}
{(1+\alpha)m-a+b-1\choose m-r-2-a}
{(1+\alpha^{-1})n+a-b-1\choose n-r-2-b}x^my^n\\[5pt]
&=\frac{uv(u-1)^{r+2}(v-1)^{r+2}} {(1+\alpha-\alpha
u)(1+\alpha^{-1}-\alpha^{-1}v)(uv-u-v)^2}.
\end{align*}
On the other hand, for $-r\leq k\leq r$, we have
\begin{align*}
\sum_{m,n=-r}^{\infty}{(1+\alpha)m\choose
m-k}{(1+\alpha^{-1})n\choose n-k}x^my^n
=\frac{uv(u-1)^{k}(v-1)^{k}} {(1+\alpha-\alpha
u)(1+\alpha^{-1}-\alpha^{-1}v)}.
\end{align*}
It's routine to verify the following identity:
\begin{align*}
G_r(x,y)+G_{-r-2}(x,y) &=\frac{2uv(u-1)(v-1)}
{(1+\alpha-\alpha u)(1+\alpha^{-1}-\alpha^{-1}v)(uv-u-v)^2}\\[5pt]
&+\sum_{k=-r}^{r}(r-|k|+1)\frac{uv(u-1)^{k}(v-1)^{k}}
{(1+\alpha-\alpha u)(1+\alpha^{-1}-\alpha^{-1}v)}.
\end{align*}
The result then follows from Theorem \ref{thm:mnalpha}.
\qed
\section{Some consequences}
\subsection{Consequences of Theorem~1}

Replacing $\alpha$, $m$, and $n$ by $q/p$, $pm$, and $qn$,
respectively, in Theorem \ref{thm:mnalpha}, we obtain
\begin{cor}\label{cor:pqmn}
For positive integers $m$, $n$, $p$, and $q$, there holds
\begin{align*}
\sum_{a=1}^{pm}\sum_{b=1}^{qn}&{pm+qm-a+b-1\choose pm-a}{pn+qn+a-b-1\choose qn-b}\\[5pt]
&=\frac{pqmn}{(p+q)(m+n)}{pm+qm\choose pm}{pn+qn\choose pn}.
\end{align*}
\end{cor}
Exchanging $p$ and $m$, and $q$ and $n$, respectively,
Corollary \ref{cor:pqmn} may be written as follows:
\begin{align*}
\sum_{a=1}^{pm}\sum_{b=1}^{qn}&{pm+qm-a+b-1\choose b-1}{pn+qn+a-b-1\choose a-1}
\nonumber\\[5pt]
&=\frac{pqmn}{(p+q)(m+n)}{pm+pn\choose pm}{qm+qn\choose qm}.
\end{align*}
By the Chu-Vandermonde formula, we have
\begin{align*}
&\sum_{a=1-pn}^{pm}\sum_{b=1}^{qn}
{pm+qm-a+b-1\choose pm-a}{pn+qn+a-b-1\choose qn-b}\\[5pt]
&=\frac{pqn}{p+q}{pm+qm+pn+qn\choose pm+pn}.
\end{align*}
Therefore, by Corollary \ref{cor:pqmn}, we have
\begin{align*}
&\sum_{a=1-pn}^{0}\sum_{b=1}^{qn}
{pm+qm-a+b-1\choose pm-a}{pn+qn+a-b-1\choose qn-b}\\[5pt]
&=\frac{pqn}{p+q}{pm+qm+pn+qn\choose pm+pn}
 -\frac{pqmn}{(p+q)(m+n)}{pm+qm\choose pm}{pn+qn\choose pn}.
\end{align*}
Replacing $a$ by $1-a$, we obtain
\begin{align*}
&\sum_{a=1}^{pn}\sum_{b=1}^{qn}
{pm+qm+a+b-2\choose pm+a-1}{pn+qn-a-b\choose qn-b}\\[5pt]
&=\frac{pqn}{p+q}{pm+qm+pn+qn\choose pm+pn}
 -\frac{pqmn}{(p+q)(m+n)}{pm+qm\choose pm}{pn+qn\choose pn}.
\end{align*}
Dividing both sides by ${pm+qm\choose pm}$, we get
\begin{align*}
&\sum_{a=1}^{pn}\sum_{b=1}^{qn}
\frac{(pm+qm+1)_{a+b-2}}{(pm+1)_{a-1}(qm+1)_{b-1}}{pn+qn-a-b\choose qn-b}\\[5pt]
&=\frac{pqn(pm+qm+1)_{pn+qn}}{(p+q)(pm+1)_{pn}(qm+1)_{qn}}
 -\frac{pqmn}{(p+q)(m+n)}{pn+qn\choose pn}.
\end{align*}
Replacing $p$, $q$, $m$, $n$, by $m$, $n$, $x$, $1$, respectively, we have
\begin{cor}For $m,n\in\mathbb{N}$, there holds
\begin{align*}
&\sum_{a=1}^{m}\sum_{b=1}^{n}
\frac{(mx+nx+1)_{a+b-2}}{(mx+1)_{a-1}(nx+1)_{b-1}}{m+n-a-b\choose m-a}\\[5pt]
&=\frac{mn(mx+nx+1)_{m+n}}{(m+n)(mx+1)_{m}(nx+1)_{n}}
 -\frac{mnx}{(m+n)(1+x)}{m+n\choose m}.
\end{align*}
\end{cor}

Letting $x\longrightarrow \infty$, we obtain
\begin{cor}For $m,n\in\mathbb{N}$, there holds
\[
\sum_{a=1}^{m}\sum_{b=1}^{n}
{m+n-a-b\choose m-a}\frac{(m+n)^{a+b-2}}{m^{a-1}n^{b-1}}\\[5pt]
=\frac{(m+n)^{m+n-1}}{m^{m-1}n^{n-1}}
 -\frac{mn}{m+n}{m+n\choose m}.
\]
\end{cor}

\subsection{Consequences of Theorem~2}
Replacing $\alpha$, $m$, and $n$ by $q/p$, $pm$, and $qn$, respectively, in
\eqref{eq:genrr-2}, we obtain
\begin{align}
&\sum_{a=1}^{pm-r-1}\sum_{b=1}^{qn-r-1}{pm+qm-a+b-1\choose pm-r-1-a}
{pn+qn+a-b-1\choose qn-r-1-b}\nonumber\\[5pt]
&+\sum_{a=0}^{pm+r}\sum_{b=0}^{qn+r}{pm+qm-a+b-1\choose pm+r-a}
{pn+qn+a-b-1\choose qn+r-b}\nonumber\\[5pt]
&=\frac{2pqmn}{(p+q)(m+n)}{pm+qm\choose pm}{pn+qn\choose pn}\nonumber\\[5pt]
&\quad\ +\sum_{k=-r}^{r}(r-|k|+1){pm+qm\choose pm-k}{pn+qn\choose qn-k}.
\label{eq:pqrsum}
\end{align}
Namely,
\begin{align}
&\sum_{a=1}^{pm-r-1}\sum_{b=1}^{qn-r-1}{pm+qm-a+b-1\choose pm-r-1-a}
{pn+qn+a-b-1\choose qn-r-1-b}\nonumber\\[5pt]
&+\sum_{a=-pm-r}^{0}\sum_{b=-qn-r}^{0}{pm+qm+a-b-1\choose pm+r+a}
{pn+qn-a+b-1\choose qn+r+b}
\nonumber\\[5pt]
&=\frac{2pqmn}{(p+q)(m+n)}{pm+qm\choose pm}{pn+qn\choose pn}\nonumber\\[5pt]
&\quad\ +\sum_{k=-r}^{r}(r-|k|+1){pm+qm\choose pm-k}{pn+qn\choose qn-k}.
\label{eq:pm-r-1}
\end{align}
By the Chu-Vandermonde formula, we have
\begin{align}
&\sum_{a=-pm-r}^{0}\sum_{b=1}^{qm-r-1}{pm+qm+a-b-1\choose pm+r+a}
{pn+qn-a+b-1\choose qn+r+b}\nonumber\\[5pt]
&+\sum_{a=-pm-r}^{0}\sum_{b=-qn-r}^{0}{pm+qm+a-b-1\choose pm+r+a}
{pn+qn-a+b-1\choose qn+r+b}
\nonumber\\[5pt]
&=\frac{(pm+r+1)q}{p+q}{(p+q)(m+n)\choose pm+pn},\label{eq:pm-r-2}
\end{align}
and
\begin{align}
&\sum_{a=-pm-r}^{0}\sum_{b=1}^{qm-r-1}{pm+qm+a-b-1\choose pm+r+a}
{pn+qn-a+b-1\choose qn+r+b}\nonumber\\[5pt]
&\sum_{a=1}^{pn-r-1}\sum_{b=1}^{qm-r-1}{pm+qm+a-b-1\choose pm+r+a}
{pn+qn-a+b-1\choose qn+r+b}\nonumber\\[5pt]
&=\frac{(qm-r-1)p}{p+q}{(p+q)(m+n)\choose pm+pn}.\label{eq:pm-r-3}
\end{align}
Summarizing \eqref{eq:pm-r-1}--\eqref{eq:pm-r-3} and replacing $r+1$ by $r$, we get
\begin{cor}
For positive integers $m$, $n$, $p$, $q$, and $r$, there holds
\begin{align}
&\sum_{a=1}^{pm-r}\sum_{b=1}^{qn-r}{pm+qm-a+b-1\choose pm-r-a}
{pn+qn+a-b-1\choose qn-r-b}\nonumber\\[5pt]
&+\sum_{a=1}^{pn-r}\sum_{b=1}^{qm-r}{pn+qn-a+b-1\choose pn-r-a}
{pm+qm+a-b-1\choose qm-r-b}\nonumber\\[5pt]
&=\frac{2pqmn}{(p+q)(m+n)}{pm+qm\choose pm}{pn+qn\choose pn}
-r{(p+q)(m+n)\choose pm+pn}\nonumber\\[5pt]
&\quad\ +\sum_{k=1-r}^{r-1}(r-|k|){pm+qm\choose pm-k}{pn+qn\choose qn-k}.
\label{eq:pm-r-4}
\end{align}
And for the case $r$ is negative, a similar formula can be deduced from
\eqref{eq:pqrsum}.
\end{cor}
 For the $p=q=1$ and $m=n=1$ cases, we obtain
the following two corollaries:
\begin{cor}For positive integers $m$, $n$, and $r$, we have
\begin{align*}
&\sum_{a=1}^{m-r}\sum_{b=1}^{n-r}{2m-a+b-1\choose m-r-a}
{2n+a-b-1\choose n-r-b}\\[5pt]
&=\frac{mn}{2(m+n)}{2m\choose m}{2n\choose n}
-\frac{r}{2}{2m+2n\choose m+n}+\frac{r}{2}{2m\choose m}{2n\choose n}
+\sum_{k=1}^{r-1}(r-k){2m\choose m-k}{2n\choose n-k}.
\end{align*}
\end{cor}

\begin{cor}For positive integers $m$, $n$, and $r$, we have
\begin{align*}
&\sum_{a=1}^{m-r}\sum_{b=1}^{n-r}{m+n-a+b-1\choose m-r-a}
{m+n+a-b-1\choose n-r-b}\\[5pt]
&=\frac{mn}{2(m+n)}{m+n\choose m}^2
-\frac{r}{2}{2m+2n\choose 2m}+\frac{r}{2}{m+n\choose m}^2
+\sum_{k=1}^{r-1}(r-k){m+n\choose m-k}{m+n\choose n-k}.
\end{align*}
\end{cor}
Furthermore, when $r=1$ and $r=2$  we obtain the following:
\begin{align*}
&\sum_{a=1}^{m-2}\sum_{b=1}^{n-2}
{2m-a+b-1\choose m-a-2}{2n+a-b-1\choose n-b-2}\\[5pt]
&={2m\choose m}{2n\choose n}+{2m\choose m-1}{2n\choose n-1}
+\frac{mn}{2(m+n)}{2m\choose m}{2n\choose n}-{2m+2n\choose m+n},
\end{align*}
\begin{align*}
&\sum_{a=1}^{m-1}\sum_{b=1}^{n-1}
{2m-a+b-1\choose m-a-1}{2n+a-b-1\choose n-b-1}\\[5pt]
&\hspace{2cm}=\frac{1}{2}{2m\choose m}{2n\choose
n}+\frac{mn}{2(m+n)}{2m\choose m}{2n\choose n}
-\frac{1}{2}{2m+2n\choose m+n}, 
\end{align*}
\begin{align*}
&\sum_{a=1}^{m-1}\sum_{b=1}^{n-1}
{m+n-a+b-1\choose m-a-1}{m+n+a-b-1\choose n-b-1}\\[5pt]
&\hspace{2cm}=\frac{1}{2}{m+n\choose
m}^2+\frac{mn}{2(m+n)}{m+n\choose m}^2 -\frac{1}{2}{2m+2n\choose 2m},
\end{align*}
and Equation~\eqref{eq:S4}.

We end this paper with one more identity of the same type:

\begin{thm}There holds
\[
\sum_{a=1}^{m}\sum_{b=1}^{n}{x+m-a+b-1\choose
n+b-1}{x+a-b-1\choose n-b} =\frac{mn}{2x+m}{2x+m\choose 2n}.
\]
\end{thm}
\pf Replacing $x$ by $-x-m+n$, one sees that the theorem is
equivalent to
\begin{equation}\label{eq:doub-xab}
2\sum_{a=1}^{m}\sum_{b=1}^{n}{x+a-1\choose n+b-1}{x+m-a\choose
n-b} =m{2x+m-1\choose 2n-1}.
\end{equation}
Now, changing $a$ to $m+1-a$ and $b$ to $1-b$, respectively, we
obtain
\[
\sum_{a=1}^{m}\sum_{b=1}^{n}{x+a-1\choose n+b-1}{x+m-a\choose n-b}
=\sum_{a=1}^{m}\sum_{b=1-n}^{0}{x+m-a\choose n-b}{x+a-1\choose
n+b-1}.
\]
So we can rewrite the left-hand side of \eqref{eq:doub-xab} as
follows:
\begin{align*}
&\sum_{a=1}^{m}\sum_{b=1}^{n}{x+a-1\choose n+b-1}{x+m-a\choose
n-b}
+\sum_{a=1}^{m}\sum_{b=1-n}^{0}{x+m-a\choose n-b}{x+a-1\choose n+b-1}\\[5pt]
&=\sum_{a=1}^{m}\sum_{b=1-n}^{n}{x+a-1\choose n+b-1}{x+m-a\choose n-b}\\[5pt]
&=m{2x+m-1\choose 2n},
\end{align*}
where the last step follows from Chu-Vandermonde's formula. \qed

\section*{Acknowledgement} We thank Christian Krattenthaler for
his helpful comments, which, in particular, led to the second
proof of Theorem~\ref{thm:mnalpha}.

\end{document}